\documentclass[12pt]{article}
\usepackage{amssymb,amsmath}
\usepackage{cases}
\usepackage{amsfonts}
\usepackage{color,xcolor}
\usepackage[left=2.0cm,right=2.0cm,top=2.0cm,bottom=2.0cm]{geometry}
\usepackage[colorlinks,citecolor=blue,urlcolor=blue]{hyperref}

\newtheorem{theorem}{Theorem}[section]
\newtheorem{corollary}{Corollary}[section]
\newtheorem{lemma}{Lemma}[section]

\newtheorem{remark}{Remark}[section]

\newcommand{\bal}{\begin{align}}
\newcommand{\bbal}{\begin{align*}}
\newcommand{\beq}{\begin{equation}}
\newcommand{\eeq}{\end{equation}}
\newcommand{\bca}{\begin{cases}}
\newcommand{\eca}{\end{cases}}
\newcommand{\pa}{\partial}
\newcommand{\fr}{\frac}
\newcommand{\na}{\nabla}
\newcommand{\De}{\Delta}

\newcommand{\lan}{\langle}
\newcommand{\ran}{\rangle}
\newcommand{\la}{\lambda}
\newcommand{\cd}{\cdot}
\newcommand{\ep}{\varepsilon}
\newcommand{\dd}{\mathrm{d}}

\newcommand{\R}{\mathbb{R}}

\newcommand{\les}{\lesssim}

\newcommand{\D}{\mathrm{div}}
\newcommand{\n}{\nabla}
\newcommand{\bi}{\Big}

\begin{document}
\title{Global large solution to the compressible Navier-Stokes equations in critical Besov space $\dot{B}^{-1}_{\infty,\infty}$}

\author{Jinlu Li$^{1}$~ Yanghai Yu$^{2}$~ Weipeng Zhu$^{3}$ and ~ Zhaoyang Yin$^{4,5}$ \thanks{E-mail: lijinlu@gnnu.cn (J. Li); yuyanghai214@sina.com (Y. Yu); mathzwp2010@163.com (W. Zhu); mcsyzy@mail.sysu.edu.cn (Z. Yin)}\\
\small $^1$\it School of Mathematics and Computer Sciences, Gannan Normal University, Ganzhou 341000, China\\
\small $^2$\it School of Mathematics and Statistics, Anhui Normal University, Wuhu, Anhui, 241002, China\\
\small $^3$\it School of Mathematics and Information Science, Guangzhou University, Guangzhou 510006, China\\
\small $^4$\it Department of Mathematics, Sun Yat-sen University, Guangzhou 510275, China\\
\small $^5$\it Faculty of Information Technology, Macau University of Science and Technology, Macau, China\\}

\date{\today}

\maketitle\noindent{\hrulefill}

{\bf Abstract:}  In this paper, we construct a class of global large solution to the compressible Navier-Stokes equations in the whole space $\R^d$. Precisely speaking, our choice of special initial data whose $\dot{B}^{-1}_{\infty,\infty}$ norm can be arbitrarily large, namely, $||u_0||_{\dot{B}^{-1}_{\infty,\infty}}\gg 1$, allows to give rise to global-in-time solution to the compressible Navier-Stokes equations.


\vskip0mm\noindent{\hrulefill}

\section{Introduction}\label{sec1}
In this paper, we focus on the compressible Navier-Stokes equations which govern the motion of a general barotropic compressible fluid in the whole space $\R^d$
\begin{eqnarray}\label{CNS}
        (\rm{CNS})~~\left\{\begin{array}{ll}
          \partial_t\varrho+\D (\varrho u)=0,\\
          \partial_t(\varrho u)+\D(\varrho u\otimes u)-\D(2\mu(\varrho)\mathbb{D}u)-\na(\lambda(\varrho)\D u)+\na P(\varrho)=0.
          \end{array}\right.
        \end{eqnarray}
Here $u=u(t,x)\in\R^d$ denotes the velocity field and $\varrho=\varrho(t,x)\in\R^+$ is the density, respectively. The density-dependent functions $\mu(\varrho)$ and $\lambda(\varrho)$ (the shear and bulk viscosity coefficients of the flow) are supposed to be smooth enough and to fulfill the standard strong parabolicity assumption:
$$\mu>0\quad\mbox{and}\quad2\mu+\lambda>0.$$
In the physical case the viscosity coefficients satisfy $2\mu+d\lambda>0$ which is a special case of the previous assumption. The strain tensor $\mathbb{D}u =\fr12(\nabla u+\nabla^{\mathsf{T}}u)$ is the symmetric part of the velocity gradient $\nabla u$. The barotropic assumption means that the pressure $P(\varrho)$
depends only upon the density $\varrho$ of fluid and the function $P$ is suitably smooth in what follows.

We supplement with initial data $(u,\varrho)|_{t=0}=(u_0,\varrho_0)$ for System \eqref{CNS} and furthermore suppose that $$\lim\limits_{|x|\to\infty}(u(t,x),\varrho(t,x))=(0,1).$$
The main difficulty in the study of the compressible fluid flows when dealing with the vacuum is that the momentum equation loses its parabolic regularizing effect. That is why in the present paper we suppose that the initial data $\varrho_0$ is a small perturbation of an equilibrium state $\bar{\varrho}$.

There is a huge literature on the studies of the compressible Navier-Stokes equations \eqref{CNS}. Here we mainly pay attention on the global existence of strong solutions to \eqref{CNS}. Kazhikhov and Shelukhin \cite{Kaz77} established the first existence result on the compressible Navier-Stokes equations in one-dimensional space for sufficiently smooth data. Serre \cite{Ser86-1,Ser86-2} and Hoff \cite{Hof87} investigated the problems for discontinuous initial data, and Mellet--Vasseur \cite{Mel2007} considered the case of dependent viscosity coefficient, where the initial density should be bounded away from the vacuum. For the multidimensional case, the first effort is due to Matsumura--Nishida \cite{Mat80} in the eighties. Actually, in three dimensions, they established the global classical solutions for initial data close to a non-vacuum equilibrium $(\bar{\varrho},0)$ under the Sobolev space setting $(\bar{\varrho}_0,u_0)\in H^3\times H^3$, and later by Hoff \cite{Hof95,Hoff95,Hof97} for discontinuous initial data. The existence of global strong solution for large initial data remains largely open for the  multidimensional compressible Navier-Stokes equations with general initial data. In the context of solutions under the assumption that small perturbations of an equilibrium state, a breakthrough was obtained by Danchin \cite{Dan00}. He addressed the global
existence issue for the isentropic compressible Navier-Stokes equations in the critical Hilbert-type Besov spaces framework. The global well-posedness result was extended to the suitable $L^p$-type Besov framework by Charve--Danchin \cite{Cha10} and Chen--Miao-Zhang \cite{Che10}, respectively. The common point between all above works is that the shear viscosity coefficient is assumed to be bounded from below by a positive constant which allows to control the gradient of the velocity, that the initial density is close to a stable constant steady state, and that the initial velocity is assumed to be small.

Although the theory is satisfactory in what concerns local time behavior and small data, many issues of global behavior of solutions for large data are far from being understood. In this research line, let us recall that there are examples of large data generating global strong solutions to the compressible Navier-Stokes equations. It should be mentioned that the remarkable work of Vaigant--Kazhikhov \cite{Vaig} who proved the
existence of global strong solution in $\R^2$ for any large initial data (with enough regularity) under the shear viscosity coefficient $\mu(\varrho)$ is a positive constant and the bulk viscosity coefficient $\lambda(\varrho)$ grows at least like $\varrho^{\beta}$ with $\beta>3$.  More
recently, B. Haspot \cite{Haspot 2017} established the existence of global strong solutions allowing for large potential part of the initial velocity in the case where $\mu(\varrho)$ depends linearly on $\varrho$ and $\lambda(\varrho)=0$. Subsequently, Danchin--Mucha \cite{Dan17} showed global existence of regular solutions to the two-dimensional compressible Navier-Stokes equations with arbitrary large initial velocity and almost constant density under the assumption of large bulk viscosity. Fang--Zhang--Zi \cite{Fan18} proved that the isentropic compressible Navier-Stokes equations admit global solutions when the initial data are close to a stable equilibrium in the sense of suitable hybrid Besov norm, where the initial velocity with arbitrary $\dot{B}^{\frac{d}{2}-1}_{2,1}$ norm of potential part and large highly oscillating are allowed. For more results of global existence of strong solutions for the multidimensional problems, we refer the reader to see \cite{Has11,Kob99,Kot12,Kot14,Lj 2019,Muc03,Muc02,Muc04,Val86} and the references therein. One can mention that Lei--Lin--Zhou \cite{Lei 2015} proved the
global well-posedness of incompressible Navier-Stokes equation with a class of large data in the energy space which includes the Beltrami flow. Motivated by this work, we continue to study the global behavior of solutions to the compressible Navier-Stokes equations \eqref{CNS} for arbitrary large initial in the present paper.
\section{Renormalized System and Statement of Main Results}\label{sec2}
Firstly, we assume that $P'(1)=1$ without loss of generality. For the sake of simplicity, we denote $$\bar{\mu}\triangleq\mu(1)\qquad \mbox{and}\qquad \bar{\lambda}\triangleq\lambda(1).$$

Our main goal is to construct the global classic solutions to \eqref{CNS} for a class of large initial data. Our results will strongly rely on the following homogeneous heat system:
\begin{eqnarray}\label{1}
        \left\{\begin{array}{ll}
          \pa_tU-\bar{\mu }\De U=0,\\
         \D U=0,\\
          U|_{t=0}=U_0.\end{array}\right.
\end{eqnarray}
As the local existence issue is nowadays well understood, we concentrate our attention on the proof of global-in-time a priori estimates under the assumption that the initial density is bounded away from zero and tends to some positive constant $\bar{\varrho}$ at infinity. We assume from now on that $\bar{\varrho}=1$ just for convenience. Denoting $\rho=\varrho-1$, hence, as long as $\varrho$ does not vanish, we can reformulate the system \eqref{CNS} equivalently as follows
\begin{eqnarray}\label{CNS-re}
        \left\{\begin{array}{ll}
          \partial_t\rho+\D u=F,\\
          \partial_tu-\bar{\mu} \De u-(\bar{\lambda}  +\bar{\mu}  )\na \D u=G,\\
          (\rho,u)|_{t=0}=(\rho_0,u_0)\end{array}\right.
        \end{eqnarray}
where
\bbal
&F=-\D(\rho u),
\\&G=-u\cd \na u-\frac{P'(\rho+1)}{\rho+1}\na \rho+\frac{2\mu'(\rho+1) }{\rho+1}\na \rho\cd\mathbb{D} u+\frac{\lambda'(\rho+1) }{\rho+1}\na\rho \D u
\\&\qquad +\Big(\frac{\mu(\rho+1) }{\rho+1}-\bar{\mu} \Big)(\De u+\na\D u)+\Big(\frac{\lambda(\rho+1) }{\rho+1}-\bar{\lambda} \Big)\na\D u.
\end{align*}

Introducing the new quantity $v=u-U$, the system \eqref{CNS-re} can be rewritten as follows
\begin{eqnarray}\label{CNS-re-1}
        \left\{\begin{array}{ll}
          \partial_t\rho+\D v=\mathbf{F},\\
          \partial_tv-\bar{\mu}  \De v-(\bar{\lambda}  +\bar{\mu}  )\na \D v+\na \rho=\mathbf{G},\\
          (\rho,v)|_{t=0}=(\rho_0,v_0)\end{array}\right.
        \end{eqnarray}
where
\bbal
&\mathbf{F}=-\rho\D v-(U+v)\cd \na \rho,
\\&\mathbf{G}=-(v+U)\cd \na (v+U)-\bi(\frac{P'(\rho+1)}{\rho+1}-1\bi)\na \rho+\frac{2\mu'(\rho+1) }{\rho+1}\na \rho\cd \mathbb{D} (v+U)
\\& \qquad+\frac{\lambda'(\rho+1) }{\rho+1}\na\rho \D v+\Big(\frac{\mu(\rho+1) }{\rho+1}-\bar{\mu} \Big)\Big(\De (v+U)+\na\D v\Big)+\Big(\frac{\lambda(\rho+1) }{\rho+1}-\bar{\lambda} \Big)\na\D v.
\end{align*}

Throughout the paper, we focus on the new system \eqref{CNS-re-1} since it is equivalent to the original system \eqref{CNS} for sufficiently smooth solutions. The main results of the paper read as follows:
\begin{theorem}\label{the1.1} Let $2\leq d\leq3$. Assume that the initial data fulfills $\rho_0=\varrho_0-1$ and $v_0=u_0-U_0$.
Let $U$ be solutions to \eqref{1} which satisfy that for some sufficiently small positive constant $\delta$ and universal constant $C$
\begin{align}\label{condition}
\bi(A_0+\int_0^\infty||U\cd \na U||_{H^3}\dd t\bi) \exp\bi\{C\int_0^\infty\bi(||U||_{W^{5,\infty}}+|| U||^2_{W^{4,\infty}}+||U\cd \na U||_{H^3}\bi)\dd t\bi\}\leq \delta,
\end{align}
where $A_0\triangleq||\rho_0,v_0||^2_{H^3}$, then the system \eqref{CNS-re-1} has a unique global solution.
\end{theorem}
Next, we establish the global classic solutions to \eqref{CNS} for a class of large initial data.

Let $a=e^{\bar{\mu} t\Delta}a_0$ be the solutions generated by the following heat equations
\bal\label{heat}
\pa_ta-\bar{\mu} \De a=0, \quad  a|_{t=0}=a_0.
\end{align}
Setting
\bbal
U=(\pa_2a, -\pa_1a)^{\mathsf{T}}\quad \mbox{for} \quad d=2 \quad \mbox{and} \quad U=
(\pa_2a, -\pa_1a, 0)^{\mathsf{T}}\quad \mbox{for} \quad d=3,
\end{align*}
then we find that $U$ solves the heat equations \eqref{1}.

Thus, Theorem \ref{the1.1} reduces to the following more exact version.
\begin{theorem}\label{the1.2} Let $2\leq d\leq3$. Assume that the initial data fulfills $\rho_0=\varrho_0-1$ and $v_0=u_0-U_0$ with
\begin{eqnarray}\label{Equ1.3}
\mathrm{supp} \ \hat{a}_0(\xi)\subset\mathcal{C}\triangleq\Big\{\xi\in\R^d: \ |\xi_1-\xi_2|\leq \ep, 1\leq |\xi|^2\leq3\Big\},\quad 0<\ep<\frac12.
\end{eqnarray}
 there exists some sufficiently small positive constant $\delta$ and universal constant $C$ such that if
\begin{align}\label{condition1}
\bi(||\rho_0,v_0||^2_{H^3}+\ep||{a}_0||_{L^2}||\hat{a}_0||_{L^1}\bi)\exp\bi\{C\bi(\ep||{a}_0||_{L^2}||\hat{a}_0||_{L^1}+||\hat{a}_0||_{L^1}+||\hat{a}_0||^2_{L^1}\bi)\bi\}\leq \delta,
\end{align}
then the system \eqref{CNS-re-1} has a unique global solution.
\end{theorem}
\begin{remark}\label{rem1.1} By choosing a special initial data, we can show that $||u_0||_{\dot{B}^{-1}_{\infty,\infty}}\gg 1$. This kind of initial data can be constructed as follows (without loss of generality, we take $v_0=0$ and $\varrho_0=1$).

For $d=3$, we set
\bbal
a_0(x_1,x_2,x_3)=\ep^{-1}\bi(\log\log\frac1\ep\bi)^{\frac12}  \chi(x_1,x_2)\phi(x_3),
\end{align*}
where the smooth functions $\chi,\phi$ satisfying $\hat{\chi}(-\xi_1,-\xi_2)=\hat{\chi}(\xi_1,\xi_2)$, $\hat{\phi}(-\xi_3)=\hat{\phi}(\xi_3)$,
\begin{align*}
\mathrm{supp} \hat{\chi}\subset \mathcal{\widetilde{C}},\quad \hat{\chi}(\xi_1,\xi_2)\in[0,1]; \quad \hat{\chi}(\xi_1,\xi_2)=1 \quad\mbox{for} \quad (\xi_1,\xi_2)\in\mathcal{\widetilde{C}}_1,
\end{align*}
and
\begin{align*}
\hat{\phi}(\xi_3)=0\quad \mbox{for} \quad|\xi_3|\in\bi[\frac{2\sqrt{2}}{3},\frac{3\sqrt{2}}{4}\bi]^c,\quad \hat{\phi}(\xi_3)\in[0,1]; \quad \hat{\phi}(\xi_3)=1 \quad\mbox{for}  \quad |\xi_3|\in\bi[\frac{\sqrt{34}}{6},\frac{\sqrt{17}}{4}\bi],
\end{align*}
where
\begin{align*}
&\mathcal{\widetilde{C}}\triangleq\Big\{\xi\in\R^2: \ |\xi_1-\xi_2|\leq \ep,\ \frac89\leq\xi^2_1+\xi^2_2\leq \frac98\Big\},
\\&\mathcal{\widetilde{C}}_1\triangleq\Big\{\xi\in\R^2: \ |\xi_1-\xi_2|\leq \frac\ep2,\ \frac{17}{18}\leq\xi^2_1+\xi^2_2\leq \frac{17}{16}\Big\} .
\end{align*}
For $d=2$, we set
\bbal
a_0(x_1,x_2)=\ep^{-1}\bi(\log\log\frac1\ep\bi)^{\frac12}  \chi(x_1,x_2),
\end{align*}
where the smooth functions $\chi$ satisfying $\hat{\chi}(-\xi_1,-\xi_2)=\hat{\chi}(\xi_1,\xi_2)$,
\begin{align*}
\mathrm{supp} \hat{\chi}\subset \mathcal{\widetilde{C}}_2,\quad \hat{\chi}(\xi)\in[0,1]\quad\mbox{and} \quad \hat{\chi}(\xi)=1 \quad\mbox{for} \quad \xi\in\mathcal{\widetilde{C}}_3,
\end{align*}
where
\begin{align*}
&\mathcal{\widetilde{C}}_2\triangleq\Big\{\xi\in\R^2: \ |\xi_1-\xi_2|\leq \ep,\ \frac{16}{9}\leq\xi^2_1+\xi^2_2\leq \frac{9}{4}\Big\},
\\&\mathcal{\widetilde{C}}_3\triangleq\Big\{\xi\in\R^2: \ |\xi_1-\xi_2|\leq \frac\ep2,\ \frac{17}{9}\leq\xi^2_1+\xi^2_2\leq \frac{17}{8}\Big\} .
\end{align*}
Then, direct calculations show that the left side of \eqref{condition} becomes
\begin{align*}
C\ep^{\frac12}\bi(\log\log \frac1\ep\bi)\exp\bi(C\log\log \frac1\ep\bi).
\end{align*}
In fact, one has
\begin{align*}
||\hat{a}_0||_{L^1}\approx \Big(\log\log\frac1\ep\Big)^\frac12\quad\mbox{and}\quad||{a}_0||_{L^2}\approx \ep^{-\fr12}\Big(\log\log\frac1\ep\Big)^\frac12.
\end{align*}
Therefore, choosing $\ep$ small enough, we deduce that the system \eqref{CNS-re-1} has a unique global solution.

Notice that $\omega_0=\pa_2u^1_0-\pa_1u^2_0=(\pa_1^2+\pa_2^2)a_0$ and $\hat{a}_0\geq 0$, we can deduce that $$\hat{w}_0=-(\xi^2_1+\xi^2_2)\hat{a}_0\leq0,$$ which implies $$||\omega_0||_{L^\infty}\approx||\hat{\omega}_0||_{L^1}.$$
Moreover, we also have $||\hat{\omega}_0||_{L^1}\gtrsim \bi(\log\log \frac1\ep\bi)^\frac12$.

Using the fact $||\omega_0||_{L^\infty}\lesssim ||u_0||_{L^\infty}$, then we have $||u_0||_{L^\infty}\gtrsim \bi(\log\log \frac1\ep\bi)^\frac12$.

Since $\mathrm{supp}\hat{u}_0\subset\bi\{|\xi|:\frac43\leq |\xi|\leq \frac32\bi\}$, then we have
$$\dot{\De}_0u_0=u_0\quad \mbox{and}\quad \dot{\De}_ju_0=0\quad\mbox{for}\quad\ j\neq 0.$$ Thus, we can conclude that
\bbal
||u_0||_{\dot{B}^{-1}_{\infty,\infty}}\approx ||u_0||_{L^\infty}\gtrsim \bi(\log\log \frac1\ep\bi)^\frac12.
\end{align*}
\end{remark}
As a by product of Theorem \ref{the1.2}, we get the global classic solutions to the incompressible Navier-Stokes equations for a class of large initial data.
\begin{corollary}\label{cor1.2} Let $2\leq d\leq3$. Assume that the initial data fulfills $v_0=u_0-U_0$. Let $U$ be solutions to \eqref{1} which satisfy \eqref{condition}, then the incompressible Navier-Stokes equations
\begin{eqnarray}\label{INS}
        (\rm{INS})~~\left\{\begin{array}{ll}
                    \partial_tu+u\cd\na u-\bar{\mu}\De u+\na P=0,\\
                   \D u=0,\\
                    u|_{t=0}=u_0.\end{array}\right.
        \end{eqnarray}
has a unique global solution.
\end{corollary}
\begin{remark} As constructed in Remark \ref{rem1.1}, the initial data with $||u_0||_{\dot{B}^{-1}_{\infty,\infty}}\gg 1$ can generate a unique global solution to the incompressible Navier-Stokes equations. This improves the classical result of Hmidi--Li \cite{Hmidi} who showed that smallness of $\dot{B}^{-1}_{\infty,\infty}$ norm of solution to d-dimensional ($d\geq3$) incompressible Navier-Stokes prevents blowups.
\end{remark}

\section{Useful Tools}\label{sec3}
\setcounter{equation}{0}
Firstly, we introduce some notations and conventions which will be used throughout this paper.
\begin{itemize}
  \item $a\lesssim b$ denotes that there are positive constants  $C$ that may vary at
different lines such that $a\leq Cb$, and we sometimes use the notation $a\approx b$ means that $a\lesssim b$ and $b\lesssim a$.
 \item We will use the abbreviated notation $||f_1,\cdots,f_n||_{X}\triangleq||f_1||_{X}+\cdots+||f_n||_{X}$ for some Banach space $X$.
  \item Let $\alpha=(\alpha_1,\cdots,\alpha_d)\in \mathbb{N}^d$ be a multi-index and $D^{\alpha}=\pa^{|\alpha|}/\pa^{\alpha_1}_{x_1}\cdots\pa^{\alpha_d}_{x_d}$ with $|\alpha|=\alpha_1+\cdots+\alpha_d$.
  \item $\lan f,g\ran$ denotes the inner product in $L^2(\R^d)$, namely, $\lan f,g\ran\triangleq\int_{\R^d}fg\dd x$.
   \item The Fourier transform of $f$ with respect to the space variable is given by
$$\hat{f}(\xi)\triangleq\int_{\R^d}e^{-ix\cd\xi}f(x)dx.$$
  \item For $m\in \mathbb{N}$, the norm of the integer order Sobolev space $H^m(\R^3)$ which is defined by
  $$||f||_{H^m(\R^d)}\triangleq\bi(\sum_{|\alpha|\leq m}||D^{\alpha}f||^2_{L^2(\R^d)}\bi)^{\fr12}=||f||_{\dot{H}^m(\R^d)}+||f||_{L^2(\R^d)}.$$
  \item The homogenous negative index Besov space $\dot{B}^{-1}_{\infty,\infty}$ equipped with the norm (see \cite{Bahouri2011})
 $$||f||_{\dot{B}^{-1}_{\infty,\infty}(\R^d)}\triangleq\sup_{j\in \mathbb{Z}}2^{-j}||\dot{\Delta}_jf||_{L^\infty(\R^d)}.$$
\end{itemize}

We give some property which will be used frequently in our proof.
\begin{lemma}\label{le2} \cite{Majda 2001} (Commutator estimates)
There hold that
\bbal
&\sum_{0<|\alpha|\leq 3}||D^{\alpha}(fg)-D^\alpha fg||_{L^2}\leq C(||f||_{{H}^{2}}||\na g||_{L^\infty}+||f||_{L^\infty}||g||_{{H}^3}),
\\&\sum_{0<|\alpha|\leq 3}||D^{\alpha}(fg)-D^\alpha fg||_{L^2}\leq C(||\na g||_{L^\infty}+||\na^3 g||_{L^\infty})||f||_{H^2}.
\end{align*}
\end{lemma}
\begin{lemma}\label{le1} \cite{Majda 2001} (Product estimates)
For $m\in \mathbb{Z}^+$ and $m\geq2$, we have
\bbal
&\sum_{|\alpha|\leq m}||D^{\alpha}(fg)||_{L^2}\leq C||f||_{H^m}||g||_{H^m},
\\&\sum_{|\alpha|\leq m}||D^{\alpha}(fg)||_{L^2}\leq C(||f||_{L^\infty}+||\na^mf||_{L^\infty})||g||_{H^m}.
\end{align*}
\end{lemma}
\begin{lemma}\label{le3}\cite{Triebel 1983} (Composition estimate)
Let $m\in \mathbb{Z}^+$ and $f\in \dot{H}^m\cap L^{\infty}$. Assume that $F\in W_{loc}^{m+2,\infty}$ with $F(0)=0$, then we have
$$
||F(f)||_{\dot{H}^m}\leq C(M)||f||_{\dot{H}^m},
$$
where the constant $C(M)$ depends on $M\triangleq\sup\limits_{k\leq m+2,|t|\leq ||f||_{L^{\infty}}}||F^{(k)}(t)||_{L^{\infty}}.$
\end{lemma}

\section{Proof of the Main Results}\label{sec4}
Before proceeding on, we introduce the following simplified notations
\bbal
A(t)\triangleq||\rho||^2_{H^3}+||v||^2_{H^3}\quad \mbox{and} \quad B(t)\triangleq||\na \rho||^2_{H^2}+||\na v||^2_{H^3}
\end{align*}
and define
\bbal
\Gamma\triangleq\sup\bi\{t\in[0,T^*):\sup_{\tau\in[0,t]}A(\tau)\leq \eta\ll1\bi\},
\end{align*}
which together with Sobolev's inequality implies that
\bal\label{y}
\sup_{\tau\in[0,t]}||\rho(\tau)||_{L^\infty}\leq C\eta\leq\frac12,
\end{align}
where $\eta$ is a small enough positive constant which will be determined later on.

Here, we emphasize that the fact which will be used often: for some smooth functions $F$ with $F(0)=0$, from \eqref{y} and Lemma \ref{le3}, there hold
$$||F(\rho)||_{L^p}\leq C||\rho||_{L^p}\quad\mbox{with}\quad p\in[1,\infty]$$
and
$$||F(\rho)||_{\dot{H}^m}\leq C||\rho||_{\dot{H}^m}\quad\mbox{with}\quad m\in\mathbb{Z}^+.$$
Now, we begin to prove the Theorem \ref{the1.1} and divide it into several parts.

{\bf Step 1: Estimation of $||\rho||_{H^3}$}.\\
Applying $D^\alpha$ to the both sides of Equs.$\eqref{CNS-re-1}_1$, taking the inner product with $D^\alpha \rho$ and summing the resulting over $|\alpha|\leq 3$ yields
\bal\label{l}
\frac12\frac{\dd}{\dd t}||\rho||^2_{H^3}+\sum_{|\alpha|\leq 3}\lan D^\alpha\D v, D^\alpha \rho \ran\triangleq I_1+I_2+I_3,
\end{align}
where
\bbal
&I_1=-\sum_{|\alpha|\leq 3}\lan D^\alpha(v\cd\na \rho), D^\alpha \rho \ran, \\
&I_2=-\sum_{0<|\alpha|\leq 3}\lan D^\alpha(U\cd\na \rho), D^\alpha \rho \ran,\\
&I_3=-\sum_{|\alpha|\leq 3}\lan D^\alpha(\rho \D v), D^\alpha \rho \ran.
\end{align*}
For the term $I_1$, we can reduce it to
\bbal
I_1=-\underbrace{\sum_{0<|\alpha|\leq 3}\lan D^\alpha(v\cd\na \rho)-v\cd \na D^{\alpha}\rho, D^\alpha \rho \ran}_{I_{1,1}}-\underbrace{\sum_{0<|\alpha|\leq 3}\lan v\cd \na D^{\alpha}\rho, D^{\alpha} \rho \ran}_{I_{1,2}}-\underbrace{\lan v\cd \na \rho, \rho \ran}_{I_{1,3}}.
\end{align*}
Notice that the embedding $H^2(\R^d)\hookrightarrow L^\infty(\R^d)$ holds for $d=2,3$. Then, by the H\"{o}lder inequality and Lemma \ref{le2}, we obtain
\bbal
&|I_{1,1}|\les\bi(||\na \rho||_{H^2}||\na v||_{L^\infty}+||v||_{H^3}||\na \rho||_{L^\infty}\bi)||\na \rho||_{H^2}\les ||v||_{H^3}||\na \rho||_{H^2}^2, \\
&|I_{1,2}|\les ||\D v||_{L^\infty}||\na \rho||^2_{H^2}\les ||v||_{H^3}||\na \rho||_{H^2}^2,
\\&|I_{1,3}|\les  ||\na v||_{L^2}||\rho||_{L^6}||\rho||_{L^3} \les ||\na v||_{L^2}||\rho||_{H^3}||\na \rho||_{L^2},
\end{align*}
where we have used the facts
\bal\label{es-23}
&||\rho||_{L^6(\R^3)}\les ||\n \rho||_{L^2(\R^3)},
\\&||\rho||_{L^6(\R^2)}\les ||\rho||^{\fr13}_{L^2(\R^2)}||\n \rho||^{\fr23}_{L^2(\R^2)}\quad\mbox{and}\quad\quad||\rho||_{L^3(\R^2)}\les ||\rho||^{\fr23}_{L^2(\R^2)}||\n \rho||^{\fr13}_{L^2(\R^2)},\nonumber
\end{align}
which implies
\bal\label{es-I1}
|I_1|\les A^{\fr12}(t)B(t).
\end{align}
For the term $I_2$, due to the fact $\D U=0$, by the H\"{o}lder inequality and Lemma \ref{le2}, we get
\bal\label{es-I2}
|I_2|=&~\bi|\sum_{0<|\alpha|\leq 3}\lan D^\alpha(U\cd\na \rho)-U\cd \na D^{\alpha}\rho, D^\alpha \rho \ran\bi|\nonumber\\
\les&~ \bi(||\na U||_{L^\infty}+||\na^3U||_{L^\infty}\bi)||\na \rho||^2_{H^2}\nonumber\\
\les&~ ||\na U,\na^3U||_{L^\infty}A(t).
\end{align}
For the term $I_3$, using the H\"{o}lder inequality, \eqref{es-23}  and Lemma \ref{le1}, we have
\bal\label{es-I3}
|I_3|\les ||\na v||_{H^3}||\rho||_{H^3}||\na \rho||_{H^2}\les A^{\fr12}(t)B(t).
\end{align}
Inserting \eqref{es-I1}--\eqref{es-I3} into \eqref{l} yields
\bal\label{lj}
\frac12\frac{\dd}{\dd t}||\rho||^2_{H^3}-\sum_{|\alpha|\leq 3}\lan D^\alpha v, D^\alpha \nabla\rho \ran\les A^{\fr12}(t)B(t)+||\na U,\na^3U||_{L^\infty}A(t).
\end{align}
{\bf Step 2: Estimation of $||v||_{H^3}$}.\\
Applying $D^\alpha$ in Equs. $\eqref{CNS-re-1}_2$ and doting the resulting equations with $D^\alpha v$, then summing them over $|\alpha|\leq 3$ gives
\bal\label{l1}
\frac12\frac{\dd}{\dd t}||v||^2_{H^3}+\bar{\mu} ||\na v||^2_{H^3}+(\bar{\lambda} +\bar{\mu} )||\D v||^2_{H^3}+\sum_{|\alpha|\leq 3}\lan D^\alpha\na \rho, D^\alpha v\ran\triangleq\sum^{12}_{i=1}J_i,
\end{align}
where
\bbal
&J_1=-\sum_{|\alpha|\leq 3}\lan D^\alpha(v\cd\na v), D^\alpha v \ran, \quad J_2=-\sum_{0<|\alpha|\leq 3}\lan D^\alpha(U\cd\na v), D^\alpha v \ran,
\\&J_3=-\sum_{|\alpha|\leq 3}\lan D^\alpha(v\cd\na U), D^\alpha v \ran, \quad J_4=-\sum_{|\alpha|\leq 3}\lan D^\alpha(U\cd \na U), D^\alpha v \ran,
\\&J_5=\sum_{|\alpha|\leq 3}\bi\lan D^\alpha\bi(\frac{2\mu'(\rho+1) }{1+\rho}\na\rho\cd \mathbb{D}v\bi), D^\alpha v \bi\ran, \quad J_6=\sum_{|\alpha|\leq 3}\bi\lan D^\alpha\bi(\frac{2\mu'(\rho+1) }{1+\rho}\na\rho\cd \mathbb{D}U\bi), D^\alpha v \bi\ran,
\\&J_7=\sum_{|\alpha|\leq 3}\bi\lan D^\alpha\bi(\frac{\lambda'(\rho+1) }{1+\rho}\na\rho \D v\bi), D^\alpha v \bi\ran, \quad J_8=-\sum_{|\alpha|\leq 3}\bi\lan D^\alpha\bi\{\bi(\frac{P'(\rho+1)}{\rho+1}-1\bi)\na \rho\bi\}, D^\alpha v \bi\ran,
\\&J_{9}=\sum_{|\alpha|\leq 3}\bi\lan D^\alpha\bi\{\bi(\frac{\mu(\rho+1) }{\rho+1}-\bar{\mu}\bi)\De v\bi\}, D^\alpha v \bi\ran, \quad J_{10}=\sum_{|\alpha|\leq 3}\bi\lan D^\alpha\bi\{\bi(\frac{\mu(\rho+1) }{\rho+1}-\bar{\mu} \bi))\De U\bi\}, D^\alpha v \bi\ran,
\\&
J_{11}=\sum_{|\alpha|\leq 3}\bi\lan D^\alpha\bi\{\bi(\frac{\mu(\rho+1) }{\rho+1}-\bar{\mu} \bi)\na\D v\bi\}, D^\alpha v \bi\ran,
J_{12}=\sum_{|\alpha|\leq 3}\bi\lan D^\alpha\bi\{\bi(\frac{\lambda(\rho+1) }{\rho+1}-\bar{\lambda} \bi)\na\D v\bi\}, D^\alpha v \bi\ran.
\end{align*}
We rewrite $J_1$ as
\bbal
J_1=-\underbrace{\sum_{0<|\alpha|\leq 3}\lan D^\alpha(v\cd\na v)-v\cd \na D^{\alpha}v, D^\alpha v \ran}_{J_{1,1}}-\underbrace{\sum_{0<|\alpha|\leq 3}\lan v\cd \na D^{\alpha}v, D^{\alpha} v \ran}_{J_{1,2}}-\underbrace{\lan v\cd \na v, v \ran}_{J_{1,3}}.
\end{align*}
Taking the similar argument as for $I_1$ yields
\bbal
&|J_{1,1}|\les\bi(||\na v||_{H^2}||\na v||_{L^\infty}+||v||_{H^3}||\na v||_{L^\infty}\bi)||\na v||_{H^2}\les ||v||_{H^3}||\na v||_{H^2}^2, \\
&|J_{1,2}|\les ||\D v||_{L^\infty}||\na v||^2_{H^2}\les ||v||_{H^3}||\na v||_{H^2}^2,\\
&|J_{1,3}|\les ||\n v||_{L^2}||v||_{L^6}||v||_{L^3}\les ||v||_{H^3}||\na v||_{L^2}^2,
\end{align*}
which implies
\bal\label{es-J1}
|J_1|\les A^{\frac12}(t)B(t).
\end{align}
Similarly, for $J_2$ and $J_3$, we have
\bal\begin{split}\label{es-J2}
&|J_2|+|J_3|\les ||\na U,\na^3U||_{L^\infty}||v||^2_{H^3}+||\na U,\na^4U||_{L^\infty}||v||^2_{H^3}\les ||\na U,\na^3U,\na^4U||_{L^\infty}A(t).
\end{split}\end{align}
To estimate $J_4$, according to H\"{o}lder's inequality, we have
\bal\label{es-J4}
|J_4|\leq ||U\cd \na U||_{H^3}||v||_{H^3}\les ||U\cd \na U||_{H^3}+||U\cd \na U||_{H^3}A(t).
\end{align}
For $J_5$, we write is as
\bbal
J_5&=-\underbrace{\sum_{0<|\alpha|\leq 3}\lan D^\alpha\bi(\frac{2\mu'(\rho+1) }{1+\rho}\na\rho\cd \mathbb{D}v\bi), D^\alpha v \ran}_{J_{5,1}}-\underbrace{\lan \bi(\frac{2\mu'(\rho+1) }{1+\rho}\na\rho\cd \mathbb{D}v\bi), v \ran}_{J_{5,2}}.
\end{align*}
By the H\"{o}lder inequality and Lemma \ref{le3}, we have
\bbal
|J_{5,1}|\les&~ \bi|\sum_{0<|\alpha|\leq 3}\lan D^{\alpha-e_i}\bi(\frac{\mu'(\rho+1) }{1+\rho}\na\rho\cd \mathbb{D}v\bi), \partial_{x_i}D^\alpha v \ran\bi|\\
\les&~\bi\|\frac{\mu'(\rho+1) }{1+\rho}\na\rho\cd \mathbb{D}v\bi\|_{{H}^2}||\nabla v||_{H^3}\\
\les&~\bi\|\frac{\mu'(\rho+1) }{1+\rho}\na\rho\bi\|_{H^2}||\mathbb{D}v||_{H^2}||\nabla v||_{H^3}\\
\les&~\bi(||\rho||_{H^3}+||\rho||^2_{H^3}\bi)||\na v||^2_{H^3}, \\
|J_{5,2}|\les&~ \bi\|\frac{\mu'(\rho+1) }{1+\rho}\bi\|_{L^\infty}||\na\rho||_{L^2} ||\n v\|_{L^2}||v||_{L^\infty}\\
\les&~\bi(\mu'(1)+\bi\|\frac{\mu'(\rho+1) }{1+\rho}-\mu'(1)\bi\|_{L^\infty}\bi)||v||_{H^3}||\na \rho||_{H^2}||\na v||_{H^3}\\
\les&~||v||_{H^3}||\na \rho||_{H^2}||\na v||_{H^3},
\end{align*}
which implies
\bal\label{es-J5}
|J_5|\les \bi(A^{\frac12}(t)+A(t)\bi)B(t).
\end{align}
For $J_6$, we write it as
\bbal
J_6&=-\underbrace{\sum_{0<|\alpha|\leq 3}\lan D^{\alpha-e_i}\bi(\frac{2\mu'(\rho+1) }{1+\rho}\na\rho\cd \mathbb{D}U\bi), D^\alpha \partial_{x_i}v \ran}_{J_{6,1}}-\underbrace{\lan \bi(\frac{2\mu'(\rho+1) }{1+\rho}\na\rho\cd \mathbb{D}U\bi), v \ran}_{J_{6,2}}.
\end{align*}
By taking advantage of the H\"{o}lder inequality and Lemmas \ref{le1} and \ref{le3}, we have
\bbal
|J_{6,1}|\les&~\bi\|\frac{\mu'(\rho+1) }{1+\rho}\na\rho\cd \mathbb{D}U\bi\|_{{H}^2}||\nabla v||_{H^3}\\
\les&~ ||\na U,\na^3U||_{L^\infty}\bi(||\rho||_{H^3}+||\rho||^2_{H^3}\bi)||\na v||_{H^3}, \\
|J_{6,2}|\les&~ ||\na U||_{L^\infty}||v||_{H^3}||\na \rho||_{H^2},
\end{align*}
which implies
\bal\label{es-J6}
|J_6|\leq  C||\na U,\na^3U||^2_{L^\infty}A(t)(1+A^2(t))+\frac{\bar{\mu}}2||\na v||^2_{H^3}.
\end{align}
Similar argument as $J_5$, we have
\bal\label{es-J7}
|J_7|\les\bi(A^{\frac12}(t)+A(t)\bi)B(t).
\end{align}
For $J_8$, we also have
\bbal
J_8=-\underbrace{\sum_{0<|\alpha|\leq 3}\lan D^\alpha\bi\{\bi(\frac{P'(\rho+1)}{\rho+1}-1\bi)\na \rho\bi\}, D^\alpha v \ran}_{J_{8,1}}-\underbrace{\lan \bi(\frac{P'(\rho+1)}{\rho+1}-1\bi)\na \rho, v \ran}_{J_{8,2}}.
\end{align*}
According to the H\"{o}lder inequality and Lemma \ref{le3}, we have
\bbal
|J_{8,1}|\les&~ ||\rho||_{H^2}||\na\rho||_{H^2}||\na v||_{H^3}, \\
\rm{For~ d=2}\quad |J_{8,2}|\les&~ \bi\|\frac{P'(\rho+1)}{\rho+1}-1\bi\|_{L^6}||\na \rho||_{L^2}||v||_{L^3}\\
\les&~||\na \rho||_{L^2}||\rho||_{L^6}||v||_{L^3}\\
\les&~ ||\rho,v||_{L^2}||\na \rho,\na v||^2_{L^2},\\
\rm{For~ d=3}\quad |J_{8,2}|\les&~\bi\|\frac{P'(\rho+1)}{\rho+1}-1\bi\|_{\dot{H}^1}||\na \rho||_{L^2}||v||_{L^3}\les ||\na \rho||^2_{H^2}||v||_{H^3},
\end{align*}
then, we can deduce that
\bal\label{es-J8}
|J_8|\leq A^{\frac12}(t)B(t).
\end{align}
For $J_{9}$, performing integrating by parts, we write it as
\bbal
J_{9}=&\underbrace{\sum_{0<|\alpha|\leq 3}\lan D^{\alpha-e_i}\bi\{\bi(\frac{\mu(\rho+1) }{\rho+1}-\bar{\mu} \bi)\cdot\Delta v\bi\}, D^\alpha \partial_{x_i} v \ran}_{J_{9,1}}
\\&-\underbrace{\sum^d_{i=1}\lan \bi\{\bi(\frac{\mu(\rho+1) }{\rho+1}-\bar{\mu} \bi)\partial_{x_i} v\bi\}, \partial_{x_i} v \ran}_{J_{9,2}}-\underbrace{\sum^d_{i=1}\lan \partial_{x_i}\bi(\frac{\mu(\rho+1) }{\rho+1}-\bar{\mu} \bi)\cdot\partial_{x_i} v,  v \ran}_{J_{9,3}}.
\end{align*}
From Lemma \ref{le3}, we have
\bbal
|J_{9,1}|\les&~ ||\rho||_{H^3}||\na v||^2_{H^3},\\
|J_{9,2}|+|J_{9,3}|\les&~\bi\|\frac{\mu(\rho+1) }{\rho+1}-\bar{\mu} \bi\|_{L^6}||\na v||_{L^3}||\na v||_{L^2}\\&+\bi\|\na\bi(\frac{\mu(\rho+1) }{\rho+1}-\bar{\mu} \bi)\bi\|_{L^2}||\na v||_{L^2}||v||_{L^\infty}\\
\les&~ ||\rho||_{H^3}||\na v||^2_{H^3}+||v||_{H^3}||\na \rho||_{L^2}||\na v||_{L^2},
\end{align*}
which implies
\bal\label{es-J9}
|J_{9}|\les A^{\frac12}(t)B(t).
\end{align}
For $J_{10}$, we infer from Lemmas \ref{le1}--\ref{le3} that
\bal\label{es-J10}
|J_{10}|\les ||\rho||_{H^3}||\na^2U,\na^5U||_{L^\infty}||v||_{H^3}\les ||\na^2U,\na^5U||_{L^\infty}A(t).
\end{align}
Similar argument as that of $J_{9}$, we have
\bal\label{es-J11}
|J_{11}|+|J_{12}|\les A^{\frac12}(t)B(t).
\end{align}
Plugging \eqref{es-J1}--\eqref{es-J11} into \eqref{l1} gives
\bal\label{lj1}\begin{split}
\frac12\frac{\dd}{\dd t}||v||^2_{H^3}+&\fr{\bar{\mu}}2 ||\na v||^2_{H^3}+(\bar{\lambda} +\bar{\mu} )||\D v||^2_{H^3}+\sum_{|\alpha|\leq 3}\lan D^\alpha\na \rho, D^\alpha v\ran\\
\les&~\bi(A^{\frac12}(t)+A(t)\bi)B(t)+\bi(\sum_{i=1}^5||\na^iU||_{L^\infty}+||U\cd \na U||_{H^3}\bi)A(t)\\
&~+||\na U,\na^3U||^2_{L^\infty}A(t)\bi(1+A^2(t)\bi)+||U\cd \na U||_{H^3}.
\end{split}\end{align}
{\bf Step 3: Estimation of $\sum_{|\alpha|\leq 2}\lan D^{\alpha}v,D^\alpha \na\rho\ran$}.\\
In order to close the above estimates, we need to bound $\sum_{|\alpha|\leq 2}\lan D^{\alpha}v,D^\alpha \na\rho\ran$. To achieve this goal, performing direct calculations gives
\bal\label{ljl}
\frac{\dd}{\dd t}\sum_{|\alpha|\leq 2}\lan D^{\alpha}v, D^\alpha \na\rho \ran+||\na \rho||^2_{H^2}-||\D v||^2_{H^2}-(\bar{\lambda}+2\bar{\mu})\sum_{|\alpha|\leq 2}\lan D^\alpha\De v, D^\alpha  \nabla\rho\ran\triangleq\sum^{15}_{i=1}K_i,
\end{align}
where
\bbal
&K_1=\sum_{|\alpha|\leq 2}\lan D^\alpha(v\cd\na \rho), D^\alpha \D v \ran,
\quad K_2=\sum_{|\alpha|\leq 2}\lan D^\alpha(U\cd\na \rho), D^\alpha \D v \ran,
\\&K_3=\sum_{|\alpha|\leq 2}\lan D^\alpha(\rho\D v), D^\alpha \D v \ran,
\quad K_4=-\sum_{|\alpha|\leq 2}\lan D^\alpha(v\cd\na v), D^\alpha \na \rho \ran,
\\&K_5=-\sum_{|\alpha|\leq 2}\lan D^\alpha(U\cd\na v), D^\alpha \na \rho  \ran,
\quad K_6=-\sum_{|\alpha|\leq 2}\lan D^\alpha(v\cd\na U), D^\alpha \na \rho  \ran,
\\&K_{7}=-\sum_{|\alpha|\leq 2}\lan D^\alpha(U\cd \na U), D^\alpha \na \rho  \ran,
\quad K_8=\sum_{|\alpha|\leq 2}\lan D^\alpha\bi(\frac{2\mu'(\rho+1) }{1+\rho}\na\rho \cd\mathbb{D}v\bi), D^\alpha \na \rho  \ran,
\\&K_9=\sum_{|\alpha|\leq 2}\lan D^\alpha\bi(\frac{2\mu'(\rho+1) }{1+\rho}\na\rho \cd\mathbb{D}U\bi), D^\alpha \na \rho \ran,\quad
K_{10}=\sum_{|\alpha|\leq 2}\lan D^\alpha\bi(\frac{\lambda'(\rho+1) }{1+\rho}\na\rho \D v\bi), D^\alpha \na \rho  \ran,
\\&K_{11}=-\sum_{|\alpha|\leq 2}\lan D^\alpha\bi\{\bi(\frac{P'(\rho+1)}{\rho+1}-1\bi)\na \rho\bi\}, D^\alpha \na \rho  \ran,
\quad K_{12}=\sum_{|\alpha|\leq 2}\lan D^\alpha\bi\{\bi(\frac{\mu(\rho+1) }{\rho+1}-\bar{\mu} \bi)\De v\bi\}, D^\alpha \na\rho \ran,\\ &K_{13}=\sum_{|\alpha|\leq 2}\int D^\alpha\bi\{\bi(\frac{\mu(\rho+1) }{\rho+1}-\bar{\mu} \bi)\De U\bi\}, D^\alpha \na\rho  \ran,\quad
K_{14}=\sum_{|\alpha|\leq 2}\lan D^\alpha\bi\{\bi(\frac{\mu(\rho+1) }{\rho+1}-\bar{\mu} \bi)\na\D v\bi\}, D^\alpha \na\rho \ran,\\
&K_{15}=\sum_{|\alpha|\leq 2}\lan D^\alpha\bi\{\bi(\frac{\lambda(\rho+1) }{\rho+1}-\bar{\lambda} \bi)\na\D v\bi\}, D^\alpha \na\rho \ran.
\end{align*}
Due to the fact that $H^2(\R^d)$ is a Banach algebra, we get
\bal\label{es-k4}\begin{split}
|K_1|+|K_3|+|K_4|\leq ||v||_{H^2}||\na \rho||_{H^2}||\na v||_{H^2}+||\rho||_{H^2}||\na v||^2_{H^2}\les A^{\frac12}(t)B(t).
\end{split}\end{align}
Invoking Lemmas \ref{le1} yields
\bal\label{es-K2}
|K_2|+|K_5|+|K_6|\les ||U,\na U,\na^2U,\na^3U||_{L^\infty}||\rho||_{H^3}||v||_{H^3}\les \sum_{i=0}^3||\na^iU||_{L^\infty} A(t),
\end{align}
Holder's inequality gives
\bal\label{es-K8}
|K_{7}|\leq ||U\cd \na U||_{H^2}||\na \rho||_{H^2}\les||U\cd \na U||_{H^2}+||U\cd \na U||_{H^2}A(t).
\end{align}
The terms $K_{8}$ and $K_{10}$ can be similarly estimated as that of $J_{5}$, we have
\bal\label{es-K7}\begin{split}
|K_8|+|K_{10}|\les&~\bi(\bi\|\frac{\mu'(\rho+1) }{1+\rho}\na\rho\cd \mathbb{D}v\bi\|_{{H}^2}+\bi\|\frac{\la'(\rho+1) }{1+\rho}\na\rho\cd \D v\bi\|_{{H}^2}\bi)||\na\rho||_{H^2}\\
&\les \bi(A^{\frac12}(t)+A(t)\bi)B(t).
\end{split}\end{align}
The terms $K_{9}$ and $K_{13}$ can be similarly estimated as that of $J_{6}$, we have
\bal\label{es-K9}\begin{split}
|K_9|+|K_{13}|\les&~\bi(\bi\|\frac{\mu'(\rho+1) }{1+\rho}\na\rho\cd \mathbb{D}U\bi\|_{{H}^2}+\bi\|\bi(\frac{\mu(\rho+1) }{\rho+1}-\bar{\mu} \bi)\De U\bi\|_{{H}^2}\bi)||\nabla \rho||_{H^2}\\
\les&~ ||\na U,\na^2 U,\na^3 U,\na^4U||_{L^\infty}\bi(||\rho||_{H^3}+||\rho||^2_{H^3}\bi)||\nabla \rho||_{H^2}\\
\les&~ \sum_{i=1}^4||\na^iU||^2_{L^\infty}A(t)(1+A^2(t))+\frac12||\nabla \rho||^2_{H^2}.
\end{split}\end{align}
From Lemmas \ref{le1}--\ref{le3} and Holder's inequality, we get
\bal\label{es-K11}
|K_{11}|+|K_{12}|+|K_{14}|+|K_{15}|\les&~ ||\nabla \rho||_{H^2}\bi\{\bi\|\bi(\frac{P'(\rho+1)}{\rho+1}-1\bi)\na \rho\bi\|_{{H}^2}+\bi\|\bi(\frac{\mu(\rho+1) }{\rho+1}-\bar{\mu} \bi)\De v\bi\|_{{H}^2}\nonumber\\
&+\bi\|\bi(\frac{\mu(\rho+1) }{\rho+1}-\bar{\mu} \bi)\na\D v\bi\|_{{H}^2}+\bi\|\bi(\frac{\lambda(\rho+1) }{\rho+1}-\bar{\lambda} \bi)\na\D v\bi\|_{{H}^2}\bi\}\nonumber\\
\les&~||\rho||_{H^3}||\na \rho||^2_{H^2}+||\rho||_{H^3}||\na \rho||_{H^2}||\Delta v||_{H^2}\nonumber\\
\les&~ A^{\frac12}(t)B(t).
\end{align}
Putting \eqref{es-K2}--\eqref{es-K11} together with \eqref{ljl1} implies
\bal\label{ljl1}\begin{split}
\frac{\dd}{\dd t}\sum_{|\alpha|\leq 2}\lan D^{\alpha}v, D^\alpha \na\rho \ran&+\fr12||\na \rho||^2_{H^2}-||\D v||^2_{H^2}-(\bar{\lambda}+2\bar{\mu})\sum_{|\alpha|\leq 2}\lan D^\alpha\De v, D^\alpha  \nabla\rho\ran\\
\les&~\bi(A^{\frac12}(t)+A(t)\bi)B(t)+\bi(||U||_{W^{3,\infty}}+||U\cd \na U||_{H^2}\bi)A(t)\\
&~+||\na U||^2_{W^{3,\infty}}A(t)\bi(1+A^2(t)\bi)+||U\cd \na U||_{H^2}.
\end{split}
\end{align}
{\bf Step 4: Closure of The A Priori Estimates}.\\
Now, we are in a position to close the above all estimates from {\bf Step 1--Step 3}.

Fundamental observations give that for some suitable positive constant $\gamma$
\bbal
||\rho||^2_{H^3}+||v||^2_{H^3}+\gamma\sum_{|\alpha|\leq 2}\lan D^{\alpha}v, D^\alpha \na\rho \ran\thickapprox A(t)
\end{align*}
and
\bbal
\frac{\gamma}{2}||\na \rho||^2_{H^2}+\frac{\bar{\mu}}{2} ||\na v||^2_{H^3}+(\bar{\lambda} +\bar{\mu} )||\D v||^2_{H^3}-\gamma||\D v||^2_{H^2}-\gamma(\bar{\lambda}+2\bar{\mu})\sum_{|\alpha|\leq 2}\lan D^\alpha\De v, D^\alpha  \nabla\rho\ran\thickapprox B(t).
\end{align*}
Multiplying the inequality \eqref{ljl1} by $\gamma$, combining \eqref{lj} and \eqref{lj1}, then integrating in time yields
\bal\label{es-con}
A(t)+\int_0^tB(s)\dd s\les&~\int_0^t\bi(A^{\frac12}(s)+A(s)\bi)B(s)\dd s+ \int_0^t\bi(||U||_{W^{5,\infty}}+||U\cd \na U||_{H^3}\bi)A(s)\dd s \nonumber\\
&~+\int_0^t||\na U||^2_{W^{3,\infty}}A(s)\bi(1+A^2(s)\bi)\dd s+\int_0^t||U\cd \na U||_{H^3}\dd s.
\end{align}
Choosing $\eta$ small enough such that the first term of right hand of \eqref{es-con} is absorbed.  Thus, for all $t\in[0,\Gamma]$, we infer from \eqref{es-con}
\bal\label{es-con1}
A(t)+\int_0^tB(s)\dd s\les&~ \int_0^t\bi(||U||_{W^{5,\infty}}+||\na U||^2_{W^{3,\infty}}+||U\cd \na U||_{H^3}\bi)A(s)\dd s \nonumber\\
&~+\int_0^t||U\cd \na U||_{H^3}\dd s.
\end{align}
Gronwall's inequality implies for all $t\in[0,\Gamma]$
\bal\label{es-con2}
A(t)\les&~\bi(A_0+\int_0^t||U\cd \na U||_{H^3}\dd s\bi) \exp\bi\{C\int_0^t\bi(||U||_{W^{5,\infty}}+||\na U||^2_{W^{3,\infty}}+||U\cd \na U||_{H^3}\bi)\dd s\bi\}\nonumber\\
\leq&~C\delta
\end{align}
provided that the condition \eqref{condition} holds.

Choosing $\eta=2C\delta$, thus we can get
\bbal
\sup_{\tau\in[0,t]}A(\tau)&\leq \fr\eta2 \quad\mbox{for}\quad t\leq \Gamma.
\end{align*}

So if $\Gamma<T^*$, due to the continuity of the solutions, we can obtain that there exists $0<\epsilon\ll1$ such that
\bbal
\sup_{\tau\in[0,t]}A(\tau)&\leq \eta \quad\mbox{for}\quad t\leq \Gamma+\epsilon<T^*,
\end{align*}
which is contradiction with the definition of $\Gamma$.

Thus, we can conclude $\Gamma=T^*$ and
\bbal
\sup_{\tau\in[0,t]}A(\tau)&\leq C<\infty \quad\mbox{for all}\quad t\in(0,T^*),
\end{align*}
which implies that $T^*=+\infty$. This completes the proof of Theorem \ref{the1.1}. $\Box$

{\bf Proof of Theorem \ref{the1.2}}\quad It follows from \eqref{heat} that $U=e^{\bar{\mu} t\Delta}U_0$ solves the following system
\begin{eqnarray}\label{2}
        \left\{\begin{array}{ll}
          \pa_tU-\bar{\mu }\De U=0,\\
         \D U=0,\\
          U|_{t=0}=U_0.\end{array}\right.
\end{eqnarray}
Notice that $\mathrm{supp}\ \hat{U}\subset \mathcal{C}$, then we have for all $m\in\mathbb{Z}^+$,
\bal\label{es-U}
||\nabla^m U||_{L^\infty}\leq C||U||_{L^\infty}\leq C||\hat{U}(t)||_{L^{1}}= Ce^{-\bar{\mu }t}\bi\||\xi|e^{-\bar{\mu }(|\xi|^2-1)t}\hat{a}_0\bi\|_{L^{1}}\leq Ce^{-\bar{\mu }t}||\hat{a}_0||_{L^{1}}.
\end{align}
To estimate $||U\cd \na U||_{H^3}$, direct calculations show that
\bbal
U\cd\na U^1&=U^1\pa_1U^1+U^2\pa_2U^1=(U^1+U^2)\pa_1U^1+U^2\pa_2(U^1+U^2)
\\&=(\pa_2-\pa_1)a\pa_1\pa_2a+\pa_1a\pa_2(\pa_1-\pa_2)a,\\
U\cd\na U^2&=U^1\pa_1U^2+U^2\pa_2U^2=(U^1+U^2)\pa_2U^2+U^1\pa_1(U^1+U^2)
\\&=(\pa_1-\pa_2)a\pa_1\pa_2a+\pa_2a\pa_1(\pa_2-\pa_1)a,\\
U\cd\na U^3&=0.
\end{align*}
By Lemma \ref{le3} and \eqref{Equ1.3}, we have
\bal\label{es-U1}
||U\cd \na U||_{H^3}\leq C||(\pa_1-\pa_2)a||_{L^\infty}||a_0||_{L^2}\leq C\ep e^{-\bar{\mu }t}||a_0||_{L^2}||\hat{a}_0||_{L^1}.
\end{align}
Thus, \eqref{Equ1.3} is ensured whenever \eqref{condition1} holds. We complete the proof of Corollary \ref{the1.2}. $\Box$

\section*{Acknowledgments} J. Li was supported by NSFC (No.11801090).

\end{document}